\documentclass[sn-mathphys-num]{sn-jnl}

\usepackage{epstopdf}

\usepackage{graphicx}%
\usepackage{multirow}%
\usepackage{amsmath,amssymb,amsfonts}%
\usepackage{amsthm}%
\usepackage{mathrsfs}%
\usepackage[title]{appendix}%
\usepackage{xcolor}%
\usepackage{textcomp}%
\usepackage{manyfoot}%
\usepackage{booktabs}%
\usepackage{algorithm}%
\usepackage{algorithmicx}%
\usepackage{algpseudocode}%
\usepackage{listings}%

\theoremstyle{thmstyleone}%
\newtheorem{thm}{Theorem}
\theoremstyle{thmstyletwo}%
\newtheorem{rem}{Remark}%

\theoremstyle{thmstylethree}%

\allowdisplaybreaks
\renewcommand{\thesection}{\arabic{section}}
\numberwithin{equation}{section}
\begin{document}

\title[Article Title]{An improved upper bound for the distribution of iterated Euler totient functions}


\author[1,2]{\fnm{Pei} \sur{Gao}}\email{gaopeihnu@163.com}
\author*[2]{\fnm{Qiyu} \sur{Yang}}\email{qiyuyangmath@163.com}

\affil*[1]{\orgdiv{Institute of Mathematics}, \orgname{Henan Academy of Sciences}, \orgaddress{\city{Zhengzhou}, \postcode{450046}, \state{Henan}, \country{China}}}

\affil[2]{\orgdiv{School of Mathematics and Statistics}, \orgname{Henan Normal University}, \orgaddress{ \city{Xinxiang}, \postcode{453007}, \state{Henan}, \country{China}}}


\abstract{Let $\phi(n)$ be the Euler totient function and $\phi_k(n)$ its $k$-fold iterate. In this note, we improve the upper bound for the number of positive $n\leqslant x$ such that $\phi_{k+1}(n)\geqslant cn$. Comparing with the upper bound which was obtained from Pollack's asymptotic formula of the summation of $\phi_{k+1}(n)$ for $n\leqslant x$, we have successfully increased the denominator exponent of the main term of the upper bound from $k$ to $k+1$.}

\keywords{Euler totient function, sieve method, iteration, counting function.}


\pacs[MSC Classification]{11N37}

\maketitle

\section{Introduction}
\subsection{Background}

For any positive integer $n$, let $\phi(n)$ be the Euler totient function given by
\begin{align}\label{varphin}
\phi(n)=n\prod_{p|n}
\left(1-\frac{1}{p}\right),
\end{align}
where $p$ runs over distinct primes dividing $n$. Let $\sigma(n)$ be the sum of divisors of $n$.

In 1929, Pillai \cite{Pillai} firstly investigated the iteration of the function $\phi(n)$. In 1944, Alaoglu and Erd\"{o}s \cite{Alaoglu} proved that
\begin{align}
\#\{n\leqslant x:\phi(\sigma(n))\geqslant cn\}=o(x),
\ \ \ \
\#\{n\leqslant x:\sigma(\phi(n))\leqslant cn\}=o(x).
\end{align}

In particular, in 2002, Luca and Pomerance \cite[Theorem 1]{Luca} obtained an asymptotic formula of the summation of $\sigma(\phi(n))/n$ for $n\leqslant x$.

In 1967, Erd\"{o}s \cite[Theorem 1]{Erdos1967} proved that for every $c<1/2$, arbitrarily small $\varepsilon>0$ and large $t$, we have
\begin{align}\label{1967k=1}
\frac{x}{\log x}(\log\log x)^t<
\#\{n\leqslant x: \phi_2(n)>cn\}
<
\frac{x}{\log x}(\log x)^{\varepsilon}
\end{align}
\begin{align}\label{1967k=2}
\#\{n\leqslant x: \phi_3(n)>cn\}
<
\frac{x}{(\log x)^2}(\log x)^{\varepsilon}
\end{align}
\begin{align}\label{1967k=3}
\#\{n\leqslant x: \phi_4(n)>cn\}
<
\frac{\alpha x}{(\log x)^2}
\end{align}
where $\alpha$ is an absolute constant independent of $c$.

Denote by $\log_k$ the $k$-fold iterated logarithm $\log\cdots \log$ ($k$-times) throughout the paper.

In 1990, Erd\"{o}s {\it et al} \cite[p191 (4.2)]{Erdos} applied Brun's method (see \cite{Halberstam}) to obtain the number of $n\leqslant x$ for which $p\nmid \phi_k(n)$. There are absolute constants $0<c_1,c_2,c_3\leqslant 1$. For any $x\geqslant x_0$, $p\leqslant (\log x)^2$, $k\leqslant c_3\log_2x$, we have
\begin{align}\label{pnmidphiknErdos}
\#\{n\leqslant x: p\nmid \phi_k(n)\}
\ll
\left\{
\begin{array}{ccccc}
x\exp\big\{-\frac{(c_2k^{-1}\log_2x)^k}{p-1}\big\}, &\hbox{if } p>\frac{(c_2k^{-1}\log_2x)^k}{c_1\log_2x}+1.\\
x(\log x)^{-c_1}, &\hbox{if } p\leqslant \frac{(c_2k^{-1}\log_2x)^k}{c_1\log_2x}+1.
\end{array}
\right.
\end{align}

They also showed that the ratio $\phi(n)/\phi_{k+1}(n)$ has a smooth, strictly increasing normal order:
(see \cite[Theorem 4.2]{Erdos}) Fix a positive integer $k$. There is a set $\mathcal{A}=\mathcal{A}(k)$ of asymptotic density 1 with the property that as $n\rightarrow \infty$ along $\mathcal{A}$, we have
\begin{align}
\phi(n)/\phi_{k+1}(n)\sim k!e^{k\gamma}(\log_3n)^k,
\end{align}
where $\gamma=0.57721\cdots$ is the usual Euler-Mascheroni constant.

In 2004, Indlekofer and K\'{a}tai \cite[Theorem 1]{Indlekofer} obtained an asymptotic formula of the number of $n\leqslant x$ such that $p\nmid \phi_k(n)$ in the range $p\in \left((\log_2x)^{k+\varepsilon}, (\log_2x)^{k+1-\varepsilon}\right)$, for $k\geqslant 2$, $\varepsilon>0$.

In 2011, according to Erd\"{o}s' result in \cite{Erdos}, Pollack \cite{Pollack2} proved that
\begin{align}\label{2011}
\sum_{n\leqslant x}\phi_{k+1}(n)\sim \frac{3}{k!e^{k\gamma}\pi^2}\cdot \frac{x^2}{(\log_3x)^k},\ \ \
\hbox{as $x\rightarrow \infty$}.
\end{align}

Naturally, we can apply Abel's identity (see \cite[Theorem 4.2]{Apostol2}) to obtain an asymptotic formula of the summation of $\phi_{k+1}(n)/n$ which deduced an initial upper bound for the number of positive $n\leqslant x$ such that $\phi_{k+1}(n)\geqslant cn$.

In 2024, Dixit and Bhattacharjee \cite[Theorem 1]{Dixit} deduced that
\begin{align}
\#\{n\leqslant x: \phi(\sigma(n))\geqslant cn\}\leqslant
\frac{\pi^2x}{6c\log_4x}+O\left(
\frac{x\log_3x}{(\log_x)^{1/\log_3x}\log_4x}\right).
\end{align}
They estimate the number of positive $n\leqslant x$ such that $\phi(\sigma(n))\geqslant cn$ by dividing into two equivalent and different cases, and improved the results by summing the the number of positive $n\leqslant x$ corresponding to these two cases. In our note, we use the same method to improve the upper bound of positive $n\leqslant x$ such that $\phi_{k+1}(n)\geqslant cn$.

\subsection{Main results}

\begin{thm}\label{thm2}
For every $c<{1}/{2}$, $0<c_1\leqslant 1$ and $1\leqslant k\leqslant {\log_2x}/{\log_4x}$
\begin{align}\label{@thm1}
\#\{n\leqslant x:\phi_{k+1}(n)\geqslant c n\}
\leqslant&
\frac{6}{(k-1)!e^{(k-1)\gamma}\pi^2}\cdot \frac{x}{c(\log_3x)^{k+1}}
+O\left(\frac{x\exp\left\{(\log_3x)^2\right\}}{(\log x)^{c_1}(\log_3x)^2}\right).
\end{align}
\end{thm}

\begin{rem}
In 2011, Pollack \cite{Pollack2} obtained an asymptotic formula (\ref{2011}) of the summation of $\phi_{k+1}(n)$ for $n\leqslant x$. By applying Abel's identity (see \cite[Theorem 4.2]{Apostol2}), for any $k\geqslant 0$ and large $x$, we have
\begin{align}
&\sum_{n\leqslant x}\frac{\phi_{k+1}(n)}{n}
=
1+\sum_{2^-< n\leqslant x}\frac{\phi_{k+1}(n)}{n}
\notag\\
=&
1+
\sum_{n\leqslant x}\phi_{k+1}(n)\frac{1}{x}
-
\sum_{n\leqslant 2^-}\phi_{k+1}(n)\frac{1}{2}
-
\int_{2^-}^{x}
\sum_{n\leqslant t}\phi_{k+1}(n)\left(\frac{1}{t}\right)'dt
\notag\\
=&
(1+o(1))
\left[
\frac{3}{k!e^{k\gamma}\pi^2}\frac{x}{(\log_3x)^k}
+
\frac{3}{k!e^{k\gamma}\pi^2}
\int_{2^-}^{x}
\frac{1}{(\log_3t)^k}dt
\right]
\notag\\
=&
\frac{6}{k!e^{k\gamma}\pi^2}\cdot \frac{x}{(\log_3x)^k}
+
O\left(\frac{x}{(\log_3x)^{k+1}\log_2x\log x}\right),\notag
\end{align}
where $\gamma=0.57721\cdots$ is the usual Euler-Mascheroni constant.

Then we can deduce an initial upper bound for the number of positive $n\leqslant x$ such that $\phi_{k+1}(n)\geqslant cn$. We have
\begin{align}\label{upperbounddirectly}
\#\{n\leqslant x: \phi_{k+1}(n)&\geqslant cn\}=\sum_{\substack{n\leqslant x\\ \phi_{k+1}(n)\geqslant cn}}1\leqslant \sum_{n\leqslant x}\frac{\phi_{k+1}(n)}{cn}
\notag\\
=&
\frac{6}{k!e^{k\gamma}\pi^2}\cdot\frac{x}{c(\log_3x)^k}
+
O\left(\frac{x}{c(\log_3x)^{k+1}\log_2x\log x}\right).
\end{align}

The goal of this note is to improve the exponent $k$ in the denominator of the main term of (\ref{upperbounddirectly}), and we ultimately succeed in increasing it to $k+1$, see (\ref{@thm1}).

In (\ref{@thm1}), taking $k=1,2,3$, we can separately obtain the upper bound for the number of $n\leqslant x$ such that $\phi_2(n)\geqslant cn$, $\phi_3(n)\geqslant cn$ and $\phi_4(n)\geqslant cn$. However, they all are poor than Erd\"{o}s' results in 1967, see \eqref{1967k=1}, \eqref{1967k=2}, \eqref{1967k=3}.

\end{rem}

\section{Proof of Theorem \ref{thm2}}

For $n>2$ and $k\geqslant 1$, we have $\phi_{k+1}(n)<n/2$, thus, in Theorem \ref{thm2}, $c<1/2$ is the best possible. Note that
\begin{align}
\phi_{k+1}(n)=\phi_k(n)\prod_{p|\phi_k(n)}\left(1-\frac{1}{p}\right).
\end{align}

Denote by $P(y):=\prod\limits_{p\leqslant y}p$, the product of primes up to $y$. If $P(y)|\phi_k(n)$, we have
\begin{align}\label{ffor1}
\phi_{k+1}(n)=&\phi_k(n)\prod_{p|\phi_k(n)}\left(1-\frac{1}{p}\right)
\leqslant \phi_k(n)\prod_{ p\leqslant y}\left(1-\frac{1}{p}\right)
<\frac{\phi_k(n)}{\log y}.
\end{align}

Here we use the fact (see \cite[Theorem 2.7(e)]{Montgomery}) that
\begin{align}
\prod_{p\leqslant y}\left(1-\frac{1}{p}\right)
\sim\frac{1}{e^{\gamma}\log y}, \ \  \ \hbox{as $y\rightarrow \infty$}
\end{align}
where $\gamma=0.57721\cdots$ is the usual Euler-Mascheroni constant.

Note that for any $c>0, \phi_{k+1}(n)<cn$ holds if $P(y)|\phi_k(n)$, $\phi_k(n)<\delta n$ and ${\delta}/{\log y}\leqslant c$. In other words, to obtain an upper bound of the number of $n\leqslant x$ such that $\phi_{k+1}(n)\geqslant cn$, it suffices to estimate the number of $n\leqslant x$ such that either $P(y)\nmid \phi_k(n)$ or $\phi_k(n)\geqslant\delta n$ under the condition ${\delta}/{\log y}\leqslant c$. Hence,
\begin{align}
\#\{n\leqslant x: \phi_{k+1}(n)\geqslant cn\}\leqslant
\#\{n\leqslant x:\phi_k(n)\geqslant \delta n\}+
\#\{n\leqslant x:P(y)\nmid \phi_k(n)\}.\notag
\end{align}

From the result of (\ref{upperbounddirectly}), we deduce that
\begin{align}\label{ffor4}
\#\{n\leqslant x:&\phi_k(n)\geqslant \delta n\}
=
\sum_{\substack{n\leqslant x\\ \phi_k(n)\geqslant \delta n}}
1
\leqslant
\sum_{n\leqslant x}
\frac{\phi_k(n)}{\delta n}
\notag\\
=&
\frac{6}{(k-1)!e^{(k-1)\gamma}\pi^2}\cdot\frac{x}{\delta(\log_3x)^{k-1}}
+
O\left(\frac{x}{\delta(\log_3x)^{k}\log_2x\log x}\right).
\end{align}

Following the upper bound in (\ref{pnmidphiknErdos}) and setting $y\leqslant \frac{(c_2k^{-1}\log_2x)^k}{c_1\log_2x}+1$, $k\leqslant c_3\log_2x$, we arrive at
\begin{align}\label{ffor5}
\#\{n\leqslant x:P(y)\nmid \phi_k(n)\}
\leqslant& \sum_{ p\leqslant y}\#\{n\leqslant x: p\nmid \phi_k(n)\}
=&
O\left(
\frac{x}{(\log x)^{c_1}}\cdot \frac{y}{\log y}
\right)
.
\end{align}

We now focus on the selection of parameters $\delta$ and $y$ to optimize the upper bound estimation. Our primary goal is to improve the initial upper bound of (\ref{upperbounddirectly}). The following four conditions (i), (ii), (iii), (iv) limit $k$ within a certain range as a precondition.

(i) We note that in the main term of (\ref{ffor4}), the larger the value of $\delta$ as the denominator, the better the upper bound. Given the constraint $\delta\leqslant c\log y$, we choose
\begin{align}
\delta = c\log y.
\end{align}

(ii) To improve the initial upper bound of (\ref{upperbounddirectly}), we need $\delta\gg\log_3 x$ in (\ref{ffor4}). Since $\delta = c\log y$, it implies
\begin{align}\label{yunder}
y\gg\log_2x.
\end{align}

(iii) To ensure that the main term of (\ref{ffor4}) is larger than the error term of (\ref{ffor5}), we have
\begin{align}
\frac{x}{(\log x)^{c_1}}\cdot \frac{y}{\log y}\ll \frac{x}{\delta (\log_3x)^{k-1}}.
\end{align}
Substituting $\delta = c\log y$ into the above inequality, we can obtain that
\begin{align}\label{yupper}
y\ll \frac{(\log x)^{c_1}}{(\log_3x)^{k - 1}}.
\end{align}
Combining with (\ref{yunder}), we can obtain an upper bound of $k$:
\begin{align}
k\ll\frac{\log_2x}{\log_4x}.
\end{align}

(iv) In (\ref{ffor5}),  we used Erd\"{o}s' result (\ref{pnmidphiknErdos}) which has a precondition:
\begin{align}
y\leqslant\frac{(c_2k^{-1}\log_2x)^k}{c_1\log_2x}+1.
\end{align}
Combining with (\ref{yunder}), we can derive the other condition of $k$,
\begin{align}
k\ll\log_2x.
\end{align}

Combining all the above conditions (i), (ii), (iii), (iv), we obtain a precondition of our main result: $1\leqslant k\leqslant \log_2x/\log_4x$.

Finally, to achieve a clear improvement of the initial upper bound (\ref{upperbounddirectly}), we choose $\delta = c(\log_3x)^2$ which increase the exponent $k$ to $k+1$ in the denominator of the main term of (\ref{upperbounddirectly}). Since $\delta = c\log y$, then $y=\exp\left((\log_3x)^2\right)$.

Choosing $\delta = c(\log_3x)^2$ and $y=\exp\left((\log_3x)^2\right)$, we can obtain that
\begin{align}
\#\{n\leqslant x: &\phi_{k+1}(n)\geqslant cn\}\leqslant
\#\{n\leqslant x:\phi_k(n)\geqslant \delta n\}+
\#\{n\leqslant x:P(y)\nmid \phi_k(n)\}
\notag\\
=&
\frac{6}{(k-1)!e^{(k-1)\gamma}\pi^2}\cdot\frac{x}{\delta(\log_3x)^{k-1}}
+
O\left(
\frac{x}{(\log x)^{c_1}}\cdot \frac{y}{\log y}
\right)
\notag\\
=&
\frac{6}{(k-1)!e^{(k-1)\gamma}\pi^2}\cdot \frac{x}{c(\log_3x)^{k+1}}
+O\left(\frac{x\exp\left\{(\log_3x)^2\right\}}{(\log x)^{c_1}(\log_3x)^2}\right).
\end{align}



\bibliographystyle{srtnumbered}
\bibliography{Eulertotient}

\end{document}